\documentclass{birkjour}
\usepackage{hyperref, enumerate}
\newtheorem{theorem}{Theorem}

\newtheorem{lemma}
{Lemma}
\newtheorem{proposition}[theorem]{Proposition}
\theoremstyle{definition}
\newtheorem{definition}
{Definition}
\theoremstyle{remark}
\newtheorem{remark}
{Remark}
\newtheorem{example}
{Example}
\begin{document}
	
	%
	%
	%
	%
	%
	\received{November 28, 2023}
	\revised{April 20, 2024}
	\accepted{April 29, 2024}

	%
	%

	\title[Clairaut Semi-invariant Riemannian Maps]{Clairaut Semi-invariant Riemannian Maps to K\"ahler Manifolds}
	
	\author[Murat Polat]{Murat Polat}
	\address{Department of Mathematics, Faculty of Science\\ Dicle University\\ Sur\\ Diyarbakir 21280\\ Turkey}
	\email{murat.polat@dicle.edu.tr}
	
	\author[Kiran Meena]{Kiran Meena}
	\address{Department of Mathematics, Lady Shri Ram College For Women\\ University of Delhi\\ Lajpat Nagar\\ New Delhi 110024\\ India}
	\email{kirankapishmeena@gmail.com}
	
	
	
	
	\begin{abstract}
		In this paper, first, we recall the notion of Clairaut Riemannian map (CRM) ${F}$ using a geodesic curve on the base manifold and give the Ricci equation. We also show that if base manifold of CRM is space form then leaves of $(ker{F}_\ast)^\perp$ become space forms and symmetric as well. Secondly, we define Clairaut semi-invariant Riemannian map (CSIRM) from a Riemannian manifold $(M, g_{M})$ to a K\"ahler manifold $(N, g_{N}, P)$ with a non-trivial example. We find necessary and sufficient conditions for a curve on the base manifold of semi-invariant Riemannian map (SIRM) to be geodesic. Further, we obtain necessary and sufficient conditions for a SIRM to be CSIRM. Moreover, we find necessary and sufficient condition for CSIRM to be harmonic and totally geodesic. In addition, we find necessary and sufficient condition for the distributions $\bar{D_1}$ and $\bar{D_2}$ of $(ker{F}_\ast)^\bot$ (which are arisen from the definition of CSIRM) to define totally geodesic foliations. Finally, we obtain necessary and sufficient conditions for $(ker{F}_\ast)^\bot$ and base manifold to be locally product manifold $\bar{D_1} \times \bar{D_2}$ and ${(range{F}_\ast)} \times {(range{F}_\ast)^\bot}$, respectively.
	\end{abstract}
	
	\maketitle
	
	\section{Introduction}
	First of all, we fix abbreviations C, RM, RMs, I, AI, SI, S, SS, HS and P for denoting Clairaut, Riemannian map, Riemannian maps, invariant, anti-invariant, semi-invariant, slant, semi-slant, hemi-slant and pointwise, respectively. The detailed geometry of Riemannian submersion was investigated in \cite{Falcitelli}. As a generalization of an isometric immersion and Riemannian submersion, Fischer introduced RM between Riemannian manifolds that satisfies the well known generalized eikonal equation $\left\Vert {F}_{\ast}\right\Vert ^{2}=rank{F} ,$ which is a bridge between geometric optics and physical optics \cite{Fischer}. Further, the geometry of RMs was investigated in \cite{Sahin1, Sahin6, Sahin_2017b, Sahin4}. 
	
	An important Clairaut's relation states that $d \sin \theta$ is constant, where $\theta$ is the angle between the velocity vector of a geodesic and a meridian, and $d$ is the distance to the axis of a surface of revolution. In 1972, Bishop defined Clairaut Riemannian submersion with connected fibers and gave a necessary and sufficient condition for a Riemannian submersion to be Clairaut Riemannian submersion \cite{Bishop}. Further, Clairaut submersions were studied in \cite{Allison}, \cite{Lee} and \cite{Meena1}. In \cite{Sahin5} and \cite{Meena}, CRMs were introduced by using geodesic curve on the total and base spaces, respectively, and obtained necessary and sufficient conditions for RMs to be CRMs. Further, \c{S}ahin gave an open problem to find characterizations for CRMs (see: \cite[p. 165, Open Problem 2]{Sahin6}).
	
	In \cite{Watson}, Watson studied almost Hermitian submersions. In \cite{Sahin3}, \c{S}ahin introduced holomorphic RMs as generalizations of holomorphic submersions and holomorphic submanifolds. Moreover, IRM, AIRM and SIRM were studied from a Riemannian manifold to a K\"ahler manifold in \cite{Akyol, Akyol1, Sahin, Sahin2}. Recently, in \cite{Yadav}, \cite{Meena_Thesis} and \cite{Meena}, Authors introduced CIRM and CAIRM from Riemannian manifolds to K\"ahler manifolds.
	
	In this paper, we study CSIRMs as generalizations of CR-submanifolds, holomorphic submersions, IRMs, CIRMs, AIRMs, CAIRMs and SIRMs from Riemannian manifolds to K\"ahler manifolds. Thus CSIRMs may have possible applications in mathematical physics. The paper is organized as follows: In Sect. \ref{sec2}, we recall some basic definitions and facts which are needed for this paper. In Sect. \ref{sec3}, first, we recall the notion of CRM using a geodesic curve on the base manifold and give the Ricci equation. We also show that if base manifold of CRM is space form then leaves of $(ker{F}_\ast)^\perp$ become space forms and symmetric as well. Secondly, we define CSIRM from a Riemannian manifold to a K\"ahler manifold. Then, we give a non-trivial example of such CSIRM. Further, we obtain a necessary and sufficient condition for a SIRM to be CSIRM. Moreover, we find necessary and sufficient condition for CSIRM to be harmonic and totally geodesic. In addition, we find necessary and sufficient condition for the distributions $\bar{D_1}$ and $\bar{D_2}$ of $(ker{F}_\ast)^\bot$ to define totally geodesic foliations. Finally, we obtain necessary and sufficient conditions for $(ker{F}_\ast)^\bot$ to be locally product manifold of $\bar{D_1}$ and $\bar{D_2}$, and for $N$ to be locally product manifold of ${(range{F}_\ast)}$ and ${(range{F}_\ast)^\bot}$. In Sect. \ref{sec4}, we give scopes of Clairaut conditions to further generalizations with almost Hermitian or K\"ahler manifolds.
	\newpage
	\section{Preliminaries}\label{sec2}
	\noindent In this section, we recall the notion of RM between Riemannian manifolds and give a brief review of basic facts.
	
	Let ${F}$ be a smooth map between Riemannian manifolds $(M, g_{M})$ and $(N,g_{N})$ such that $0 < rank{F} < \min \left\{\dim(M), \dim (N)\right\}$. We denote the kernel space of ${F}_{\ast}$ by $V_{q}=ker {F}_{\ast {q}}$ at $q \in M$ and consider the orthogonal complementary space $H_{q}= (ker {F}_{\ast {q}})^{\bot}$ to $ker {F}_{\ast {q}}$ in $T_{q} M$. Then, the tangent space $T_{q} M$ of $M$ at $q$ has the decomposition
	\begin{equation*}
		T_{q} M = (ker {F}_{\ast {q}}) \oplus (ker {F}_{\ast {q}})^{\bot} = V_{q} \oplus H_{q}.
	\end{equation*}
	We denote the range of ${F}_{\ast}$ by $range{F}_{\ast}$ at $q \in M$ and consider the orthogonal complementary space $(range {F}_{\ast {q}})^{\bot}$ to $range {F}_{\ast{q}}$ in the tangent space $T_{{F} (q)}N$ of $N$ at ${F} (q) \in N$. Since $rank {F} < \min \left\{\dim (M), \dim (N)\right\}$, we have $(range {F}_{\ast {q}})^{\bot} \neq \left\{ 0 \right\}$. Thus, the tangent space $T_{{F} (q)} N$ of $N$ at ${F} (q) \in N$ has the decomposition
	\begin{equation*}
		T_{{F} (q)}N=(range{F}_{\ast _{q}})\oplus (range{F}_{\ast_{q}})^{\bot}.
	\end{equation*}
	Now, ${F}$ is called RM at $q \in M$ if the horizontal restriction ${F}_{\ast {q}}^{h}: (ker {F}_{\ast {q}})^{\bot} \to (range {F}_{\ast {q}})$ is a linear isometry between the spaces $((ker {F}_{\ast {q}})^{\bot}, g_{M_{{q}}}|_{(ker {F}_{\ast {q}})^{\bot}})$ and $(range {F}_{\ast {q}}, g_{N {F}(q)}|_{range {F}_{\ast {q}}})$. Equivalently,
	\begin{equation}\label{2.1}
		g_{N}({F}_{\ast}X, {F}_{\ast}Y) = g_{M}(X,Y) 
	\end{equation}
	for all $X, Y$ vector field tangent to $\Gamma (ker {F}_{\ast {q}})^{\bot}$. It follows that isometric immersions and Riemannian submersions are particular RMs with $ker {F}_{\ast} = \left\{ 0 \right\}$ and $(range {F}_{\ast})^{\bot} = \left\{ 0 \right\}$, respectively. The differential map ${F}_{\ast}$ of ${F}$ can be viewed as a section of bundle $\hom (TM, {F}^{-1} TN) \to M$, where ${F}^{-1} TN$ is the pullback bundle whose fibers at $q \in M$ is $\left({F}^{-1} TN\right)_{q} = T_{{F} (q)}N$, $q \in M$. The bundle $\hom (TM, {F}^{-1} TN)$ has a connection $\nabla$ induced from the Levi-Civita connection $\nabla^{M}$ and the pullback connection $\overset{N}{\nabla^{{F}}}$. Then, the second fundamental form of ${F}$ is defined as \cite{Nore} 
	\begin{equation}\label{2.2}
		(\nabla {F}_{\ast})(X, Y)= \overset{N}{\nabla_{X}^{{F} }} {F}_{\ast} (Y) - {F}_{\ast}(\nabla_{X}^{M}Y)
	\end{equation}
	for all $X, Y \in \Gamma (TM)$, where $\overset{N}{\nabla_{X}^{{F} }} {F}_{\ast} Y \circ {F} = \overset{N}{\nabla_{{F}_{\ast}X}^{{F}}} {F}_{\ast}Y$. In \cite{Sahin4} \c{S}ahin showed that $(\nabla {F}_{\ast})(X, Y) = (\nabla {F}_{\ast})(Y, X)$. In \cite{Sahin}, \c{S}ahin proved for all $X, Y \in \Gamma (ker {F}_{\ast})^{\bot}$
	\begin{equation}\label{2.3}
		(\nabla {F}_{\ast}) (X, Y) \in \Gamma (range{F}_{\ast})^{\bot}. 
	\end{equation}
	
	\begin{lemma}\label{lemmaforumbilicity}
		\cite{Sahin2} Let ${F}$ be a RM between Riemannian manifolds $(M, g_{M})$ and $(N,g_{N})$. Then, ${F} $ is umbilical RM if and only if
		\begin{equation}\label{2.4}
			(\nabla {F}_{\ast})(X,Y) = g_{M}(X, Y) H
		\end{equation}
		for all $X, Y \in \Gamma (ker {F}_{\ast})^{\bot}$ and $H$ is the mean curvature vector field of $range{F}_{\ast}$.
	\end{lemma}
	For any vector field $X$ on $M$ and any section $V$ of $(range{F}_{\ast})^{\bot}$, we have $\nabla_{X}^{{F} \bot}V$, which is the orthogonal projection of $\nabla_{X}^{N}V$ on $(range{F}_{\ast})^{\bot}$, where $\nabla^{{F} \bot}$ is linear connection on $(range{F}_{\ast})^{\bot}$ such that $\nabla^{{F} \bot} g_{N} = 0$. In addition, $(range{F}_{\ast})^{\bot}$ is totally geodesic if and only if $\nabla_{U}^{N} V = \nabla_{U}^{{F} \bot}V$ for all $U, V \in \Gamma (range{F}_{\ast})^{\bot}$ \cite{Yadav_PMD}. Now, for a RM ${F}$, we have $S_{V}$ as \cite[p. 188]{Sahin4}
	\begin{equation}\label{2.5}
		\nabla_{{F}_{\ast}X}^{N}V = - S_{V} {F}_{\ast} X + \nabla_{X}^{{F} \bot}V, 
	\end{equation}
	where $\nabla^{N}$ is Levi-Civita connection on $N$, $S_{V} {F}_{\ast} X$ is the tangential component (a vector field along ${F}$) of $\nabla_{{F}_{\ast}X}^{N} V$. Thus at $q \in M$, we have $\nabla_{{F}_{\ast} X}^{N}V(q) \in T_{{F} (q)} N, S_{V} {F}_{\ast} X \in {F}_{\ast q}(T_{q} M)$ and $\nabla_{X}^{{F} \bot} V \in ({F}_{\ast q}(T_{q}M))^{\bot}$. It is easy to see that $S_{V} {F}_{\ast} X$ is bilinear in $V$, and ${F}_\ast X$ at $q$ depends only on $V_{q}$ and ${F}_{\ast q} X_{q}$. 
	
	Let $(N, g_{N})$ be an almost Hermitian manifold \cite{Yano}, then $N$ admits a tensor $P$ of type $(1, 1)$ on $N$ such that $P^{2} = - I$ and
	\begin{equation}\label{2.8}
		g_{N}(P\tilde{X}, P\tilde{Y}) = g_{N}(\tilde{X}, \tilde{Y})
	\end{equation}
	for all $\tilde{X}, \tilde{Y} \in \Gamma (TN)$. An almost Hermitian manifold $N$ is called K\"ahler manifold if
	\begin{equation*}
		(\nabla_{\tilde{X}}^{N}P)\tilde{Y} = 0,
	\end{equation*}
	where $\nabla^{N}$ is the Levi-Civita connection on $N$.
	\begin{definition}
		\cite{Sahin2} Let ${F}$ be a RM from a Riemannian manifold $(M, g_{M})$ to an almost Hermitian manifold $(N, g_{N}, P)$. Then we say that ${F}$ is a SIRM at $q \in M$ if there are subbundles $D_{1}$ and $D_{2}$ in $range{F}_{\ast}$ such that $P(D_{1}) = D_{1}~ \& ~P(D_{2}) \subseteq (range{F}_{\ast})^{\bot}$. If ${F}$ is a SIRM at every point $q \in M$, then we say that ${F}$ is a SIRM.
	\end{definition}
	\noindent Then, for ${F}_\ast X \in \Gamma (range{F}_{\ast})$, we write
	\begin{equation}\label{3.3}
		P {F}_{\ast} X = \phi {F}_{\ast}X + \omega {F}_{\ast}X,
	\end{equation}
	where $\phi {F}_{\ast} X \in \Gamma (D_{1})$ and $\omega {F}_{\ast} X \in \Gamma (PD_{2})$. Also for ${F}_{\ast} X \in \Gamma(D_1)$ and ${F}_\ast Y \in \Gamma(D_2)$, we have $g_{N} ({F}_\ast X, {F}_{\ast} Y) = 0$. Thus, we have two orthogonal distributions $\bar{D_1}$ and $\bar{D_2}$ such that
	\[(ker {F}_\ast)^\bot = \bar{D_1} \oplus \bar{D_2}.\] In addition, for $V \in \Gamma((range {F}_{\ast})^{\bot})$, we write
	\begin{equation}\label{3.4}
		PV = BV + CV, 
	\end{equation}
	where $BV \in \Gamma (D_{1})$ and $CV \in \Gamma (\eta)$. Here $\eta$ is the complementary orthogonal distribution to $\omega (D_{2})$ in $(range {F}_{\ast})^{\bot}$. It is easy to see that $\eta$ is invariant with respect to $P$.
	
	\section{CSIRMs from Riemannian Manifolds to K\"ahler Manifolds}\label{sec3}
	
	\noindent In this section, we introduce CSIRM from a Riemannian manifold to a K\"ahler manifold and investigate the geometry.
	
	The notion of CRM comes from a geodesic curve on a surface of revolution. \c{S}ahin defined the CRM using a geodesic curve on the total manifold in \cite{Sahin5}. Recently, CRM was defined by Meena and Yadav using a geodesic curve on the base manifold in \cite{Meena} as:
	
	\begin{definition}\label{DEF1}
		\cite{Meena} A RM ${F}$ between Riemannian manifolds $(M, g_{M})$ and $(N,g_{N})$ is called CRM if there is a function $s:N\to \mathbb{R}^{+}$ such that for every geodesic $\beta$ on $N$, the function $(s\circ \beta )\sin \theta (t)$ is constant, where, ${F}_{\ast} X \in \Gamma (range{F}_{\ast})$ and $V\in \Gamma (range{F}_{\ast})^{\bot}$ are components of $\dot{\beta}(t)$, and $\theta (t)$ is the angle between $\dot{\beta}(t)$ and $V$ for all $t$.
	\end{definition}
	\noindent The authors Meena and Yadav obtained the following necessary and sufficient conditions for a RM to be CRM:
	\begin{theorem}\label{TH1}
		\cite{Meena} Let ${F}$ be a RM  between Riemannian manifolds $(M, g_{M})$ and $(N,g_{N})$ such that $(range{F}_{\ast})^{\bot}$ is totally geodesic and leaves of $range{F}_\ast$ are connected, and let $\alpha, \beta = {F} \circ \alpha$ be geodesic curves on $M$ and $N$, respectively. Then ${F}$ is CRM with $s = e^{f}$ if and only if any one of the following conditions holds:
		\begin{enumerate}[$(i)$]
			\item $S_{V} {F}_{\ast} X = -V(f) {F}_{\ast}X$, where ${F}_{\ast} X \in \Gamma(range{F}_{\ast})$ and $V \in \Gamma (range{F}_{\ast})^{\bot}$ are components of $\dot{\beta}(t)$, i.e. $\dot{\beta}(t) = {F}_{\ast} X + V$.
			
			\item ${F}$ is umbilical map, and has $H = - \nabla^{N}f$, where $f$ is a smooth function on $N$ and $H$ is the mean curvature vector field of $range{F}_{\ast}$.
		\end{enumerate}
	\end{theorem}
	
	\begin{remark}
		Note that from now onward, we are assuming the map ${F}$ such that $(range {F}_{\ast})^{\bot}$ is totally geodesic. This is needed to investigate the geometry of CRMs.
	\end{remark}
	
	\begin{theorem}
		Let ${F}$ be a CRM with $s=e^{f}$ between Riemannian manifolds $(M, g_{M})$ and $(N, g_{N})$. Then the Ricci equation is
		\begin{equation*}
			[S_W, S_V] F_\ast X = 0,
		\end{equation*}
		where $X, Y \in \Gamma(kerF_\ast)^\perp$ and $W, V \in \Gamma (range{F}_{\ast})^{\bot}$.
	\end{theorem}
	
	\begin{proof}
		For a RM, we have \cite{Sahin4}
		\begin{align*}
			g_N(R^N(F_\ast X, F_\ast Y)V, W)& = g_N(R^{F \perp}(F_\ast X, F_\ast Y)V, W) \\&-g_N((\nabla F_\ast)(X, {}^\ast F_\ast(S_V F_\ast Y)), W) \\& + g_N((\nabla F_\ast)(Y, {}^\ast F_\ast(S_V F_\ast X)), W),
		\end{align*}
		where $X, Y \in \Gamma(kerF_\ast)^\perp$ and $W, V \in \Gamma (range{F}_{\ast})^{\bot}$. Since $F$ is a CRM, above equation can be written as
		\begin{align*}
			g_N(R^N(F_\ast X, F_\ast Y)V, W)& = g_N(R^{F \perp}(F_\ast X, F_\ast Y)V, W) \\&+ g_M(X, {}^\ast F_\ast(S_V F_\ast Y))g_N(\nabla^N f, W)\\& - g_M(Y, {}^\ast F_\ast(S_V F_\ast X))g_N(\nabla^N f, W),
		\end{align*}
		which is equal to
		\begin{align}\label{RFperp}
			g_N(R^N(F_\ast X, F_\ast Y)V, W)& = g_N(R^{F \perp}(F_\ast X, F_\ast Y)V, W) \nonumber \\&+ g_N(F_\ast X, S_V F_\ast Y)g_N(\nabla^N f, W) \nonumber\\&- g_N(F_\ast Y, S_V F_\ast X)g_N(\nabla^N f, W).
		\end{align}
		Using Theorem \ref{TH1} in (\ref{RFperp}), we get
		\begin{align}\label{RFperpreduced}
			g_N(R^N(F_\ast X, F_\ast Y)V, W) = g_N(R^{F \perp}(F_\ast X, F_\ast Y)V, W).
		\end{align}
		We know that the Ricci equation for a RM is given by \cite{Sahin4}
		\begin{align}\label{riccieqn}
			g_N(R^N(F_\ast X, F_\ast Y)V, W)& = g_N(R^{F \perp}(F_\ast X, F_\ast Y)V, W) \nonumber\\& + g_N([S_W, S_V]F_\ast X, F_\ast Y).
		\end{align}
		Thus by (\ref{RFperpreduced}) and (\ref{riccieqn}), we obtain
		\begin{equation*}
			[S_W, S_V] F_\ast X = 0.
		\end{equation*}
	\end{proof}
	
	\begin{remark}
		A Riemannian manifold is {\it locally symmetric} if and only if $\nabla R \equiv 0$ \cite{Petersen_2016}. In particular, the space of constant sectional curvature is a locally symmetric space.
	\end{remark}
	
	\begin{theorem}
		Let ${F}$ be a CRM with $s=e^{f}$ from a Riemannian manifolds $(M, g_{M})$ to a space form $(N(k), g_{N})$ such that $\|\nabla^N f\| = 1$. If $(ker F_\ast)^\perp$ is integrable distribution then any leaf of $(ker F_\ast)^\perp$ is also space form and symmetric as well. 
	\end{theorem}
	
	\begin{proof} For a RM we have the Gauss equation \cite{Sahin4} 
		\begin{align*}
			g_M(R^M(X, Y)Z, T)& = g_N(R^N(F_\ast X, F_\ast Y)F_\ast Z, F_\ast T) \\&- g_N((\nabla F_\ast)(X, Z), (\nabla F_\ast)(Y, T)) \\&+ g_N((\nabla F_\ast)(Y, Z), (\nabla F_\ast)(X, T)),
		\end{align*}
		where $X, Y, Z ,T \in \Gamma(ker F_\ast)^\perp$. Since $N(k)$ is space form and $F$ is CRM, above equation can be written as
		\begin{align*}
			g_M(R^M(X, Y)Z, T)& = kg_N(F_\ast Y, F_\ast Z) g_N(F_\ast X, F_\ast T) \\&- kg_N(F_\ast X, F_\ast Z) g_N(F_\ast Y, F_\ast T)\\& - g_M(X, Z) g_M(Y,T) g_N(\nabla^N f, \nabla^N f) \\&+ g_M(Y, Z) g_M(X, T) g_N(\nabla^N f, \nabla^N f).
		\end{align*}
		Since we have $\|\nabla^N f\| = 1$, we write
		\begin{align*}
			g_M(R^M(X, Y)Z, T)& = (k + 1)\{g_M(Y, Z) g_M(X, T) - g_M(X, Z) g_M(Y, T)\}.
		\end{align*}
		This implies the proof.
	\end{proof}
	
	\begin{definition}
		A SIRM from a Riemannian manifold to a K\"ahler manifold is called CSIRM if it satisfies the Definition \ref{DEF1} of CRM.
	\end{definition}
	
	\begin{example}
		Let $M$ be an Euclidean space given by
		\begin{equation*}
			M = \left\{ (u_{1}, u_{2}, u_{3}, u_{4}) \in \mathbb{R}^{4} : u_{i} \neq 0~ \forall~ 1\leq i\leq 4 \right\}.
		\end{equation*}
		We describe the Riemannian metric $g_{M}$ on $M$ given by
		\begin{equation*}
			g_{M} = e^{2u_1}du_{1}^{2} + e^{2u_1} du_{2}^{2} + du_{3}^{2} + du_{4}^{2}.
		\end{equation*}
		Let $N = \left\{(v_{1}, v_{2}, v_{3}, v_{4}) \in \mathbb{R}^{4} \right\}$ be a Riemannian manifold with Riemannian metric $g_{N}$ on $N$ given by 
		\begin{equation*}
			g_{N} = e^{2v_1} dv_{1}^{2} + e^{2v_1} dv_{2}^{2} +  dv_{3}^{2} +  dv_{4}^{2}.
		\end{equation*}
		We take the complex structure $P$ on $N$ as $P(a, b, c, d) = (-b, a, -d, c)$. 
		Then, a basis of $T_q M$ is
		\begin{equation*}
			\left\{ e_{i} = \frac{\partial}{\partial u_{i}} ~\text{for $1 \leq i \leq 4$}\right\},
		\end{equation*}
		and a $P$-basis on $T_{{F}(q)} N$ is
		\begin{equation*}
			\left\{e_{j}^{\ast} = \frac{\partial}{\partial v_{j}} ~\text{for $1 \leq i \leq 4$}\right\},
		\end{equation*}
		for all $q \in M$. Now, we define a map ${F}$ between $(M, g_{M})$ and $(N, g_{N}, P)$ by
		\begin{equation*}
			{F}(u_{1}, u_{2}, u_{3}, u_{4}) = \left(u_1, u_{2}, u_{3},0\right).
		\end{equation*}
		Then, we have
		\begin{equation*}
			ker {F}_{\ast} = Span \left\{U = e_4\right\},
		\end{equation*}
		\begin{equation*}
			(ker {F}_{\ast})^{\bot} = Span\left\{X_{1} = e_1, X_{2} = e_{2}, X_{3} = e_{3}\right\}.
		\end{equation*}
		It is easy to see that ${F}_{\ast} X_{1} = e_{1}^{\ast}, {F}_{\ast} X_{2} = e_{2}^{\ast}, {F}_{\ast} X_{3} = e_{3}^{\ast}$ and $g_{M} (X_{i}, X_{j}) = g_{N} ({F}_{\ast} \left( X_{i}\right), {F}_{\ast} \left( X_{j} \right))$ for $i, j = 1, 2, 3$. Thus ${F}$ is a RM with 
		\begin{equation*}
			range {F}_{\ast} = Span \left\{e_{1}^{\ast}, e_{2}^{\ast}, e_{3}^{\ast} \right\}, (range{F}_{\ast})^{\bot} = Span \left\{e_{4}^{\ast}\right\},
		\end{equation*}
		\begin{equation*}
			D_1 = Span \left\{e_{1}^{\ast}, e_{2}^{\ast}\right\}, D_2 = Span \left\{e_{3}^{\ast}\right\}.
		\end{equation*}
		Moreover, $P {F}_{\ast} X_{1} = e_{2}^{\ast}, P {F}_{\ast} X_{2} = -e_{1}^{\ast}, P {F}_{\ast} X_{3} =  e_{4}^{\ast}$. Thus ${F}$ is a SIRM.  Since we have the only non-zero Christoffel symbols $\Gamma^{2}_{12}=\Gamma^2_{21}=1$, $\Gamma^{1}_{11} = 1$,
		$\Gamma^1_{22}= -1$ for $g_{N}$, $\nabla_{e_{4}^\ast}^N e_4^\ast =0$. This implies $(rangeF_\ast)^\perp$ is totally geodesic. Also, we can check easily that $$(\nabla_{\tilde{X}}^{N}P)\tilde{Y} = 0$$ for $\tilde{X}=a_1e_1^\ast + a_2 e_2^\ast+a_3e_3^\ast + a_4 e_4^\ast$ and $\tilde{Y}= b_1e_1^\ast + b_2 e_2^\ast+b_3e_3^\ast + b_4 e_4^\ast$. Here, $a_i, b_i\in \mathbb{R}$ for $1\leq i\leq 4$.
		
		\noindent Now for $X \in \Gamma(kerF_\ast)^\perp$, we have
		\begin{align*}
			g_M(X, X)& = g_M(\lambda_1 e_1 + \lambda_2 e_2 + \lambda_3 e_3, \lambda_1 e_1 + \lambda_2 e_2 + \lambda_3 e_3) \\&= (\lambda_1^2 + \lambda_2^2) e^{2u_1}+ \lambda_3^2
		\end{align*}
		for some $\lambda_1, \lambda_2, \lambda_3 \in \mathbb{R}$. We know that for $X \in \Gamma(kerF_\ast)^\perp$ the second fundamental form of a RM, $(\nabla {F}_{\ast})(X, X) \in \Gamma(rangeF_\ast)^\perp$, therefore we can write $$(\nabla {F}_{\ast})(X, X) = \lambda e_4^\ast$$ for some $\lambda \in \mathbb{R}$. Now it is easy to verify that $(\nabla {F}_{\ast})(X, X) = - g_{M} (X, X) \nabla^{N} f$, for $f = \frac{-\lambda v_4}{(\lambda_1^2 + \lambda_2^2) e^{2u_1}+\lambda_3^2}$; where at least one of $\lambda_1, \lambda_2, \lambda_3$ is non-zero. Hence, ${F}$ is a CSIRM from a Riemannian manifold $(M, g_{M})$ to a K\"ahler manifold $(N, g_{N}, P)$.
	\end{example}

	\begin{theorem}
		Let ${F}$ be a SIRM from a Riemannian manifold $(M, g_{M})$ to a K\"ahler manifold $(N, g_{N}, P)$ and $\alpha: I \to M$ be a geodesic curve on $M$. Then the curve $\beta = {F} \circ \alpha$ is geodesic on $N$ if and only if
		\begin{align}\label{3.5}
			&-S_{\omega {F}_{\ast} X} {F}_{\ast} X + {F}_{\ast} (\nabla_{X}^{M} {}^{\ast} {F}_{\ast} BV) - S_{CV} {F}_{\ast} X \nonumber\\& + \nabla_{V}^{N} \phi {F}_{\ast} X + \nabla_{V}^{N} BV + {F}_\ast(\nabla_X^{M} {}^{\ast}{F}_\ast \phi {F}_\ast X)= 0
		\end{align}
		and
		\begin{align}\label{3.6}
			&(\nabla {F}_\ast) (X, {}^{\ast}{F}_\ast \phi {F}_\ast X) + \nabla_{X}^{{F} \bot} \omega {F}_{\ast} X + (\nabla {F}_{\ast}) (X,{}^{\ast}{F}_{\ast} BV) \nonumber\\&+ \nabla_{X}^{{F} \bot} CV + \nabla_{V}^{{F} \bot} \omega {F}_{\ast} X + \nabla_{V}^{{F} \bot} CV = 0, 
		\end{align}
		where ${F}_{\ast} X \in \Gamma (range{F}_{\ast}), V \in \Gamma (range {F}_{\ast})^{\bot}$ are components of $\dot{\beta}(t)$ and ${}^{\ast} {F}_{\ast}$ is the adjoint map of ${F}_{\ast}$. In addition, $\nabla^{N}$ is the Levi-Civita connection on $N$ and $\nabla^{{F} \bot}$ is a linear connection on $(range{F}_{\ast})^{\bot}$.
	\end{theorem}
	
	\begin{proof}
		Let $\alpha : I \to M$ be a geodesic on $M$ and let $\beta = {F} \circ \alpha$ be a regular curve on $N$ with ${F}_{\ast} X \in \Gamma (range {F}_{\ast}), V \in \Gamma(range{F}_{\ast})^{\bot}$ are components of $\dot{\beta}(t)$. Since $N$ is K\"ahler manifold, $\nabla_{\dot{\beta}}^{N} \dot{\beta} = - P \nabla_{\dot{\beta}}^{N} P\dot{\beta}$. Thus
		\begin{equation*}
			\nabla_{\dot{\beta}}^{N} \dot{\beta} = - P \nabla_{\dot{\beta}}^{N} P \dot{\beta} = - P \nabla_{{F}_{\ast} X + V}^{N} P ({F}_{\ast} X + V),
		\end{equation*}
		which implies
		\begin{equation}\label{3.7}
			\nabla_{\dot{\beta}}^{N} \dot{\beta} = - P \left(\nabla_{{F}_{\ast} X}^{N} P {F}_{\ast} X + \nabla_{{F}_{\ast} X}^{N} P V + \nabla_{V}^{N} P {F}_{\ast} X + \nabla_{V}^{N} PV \right). 
		\end{equation}
		Using (\ref{2.5}), (\ref{3.3}) and (\ref{3.4}) in (\ref{3.7}), we have
		\begin{eqnarray}\label{3.8}
			\nabla_{\dot{\beta}}^{N} \dot{\beta} =& - P \left( \begin{array}{c} \nonumber \nabla_{{F}_{\ast} X }^{N} \phi {F}_{\ast} X + \nabla_{{F}_{\ast} X }^{N} \omega {F}_{\ast} X + \nabla_{{F}_{\ast} X }^{N} B V + \nabla_{{F}_{\ast} X }^{N} CV \\ + \nabla_{V}^{N} \phi {F}_{\ast} X + \nabla_{V}^{N} \omega {F}_{\ast} X + \nabla_{V}^{N} B V + \nabla_{V}^{N} C V
			\end{array}
			\right)\\
			= & - P \left( 
			\begin{array}{c}
				{F}_\ast(\nabla_X^{M} {}^{\ast}{F}_\ast \phi {F}_\ast X) + (\nabla {F}_\ast) (X, {}^{\ast}{F}_\ast \phi {F}_\ast X) + \nabla_{X}^{{F} \bot} C V \\- S_{\omega{F}_{\ast}X} {F}_{\ast} X + \nabla_{X}^{{F}\bot} \omega {F}_{\ast} X + \nabla_{{F}_{\ast} X}^{N} B V - S_{CV} {F}_{\ast} X \\+ \nabla_{V}^{N} \phi {F}_{\ast} X + \nabla_{V}^{N} \omega {F}_{\ast} X + \nabla_{V}^{N} B V + \nabla_{V}^{N} C V
			\end{array}
			\right).
		\end{eqnarray}
		Since $\nabla^{N}$ is Levi-Civita connection on $N$, $g_{N}(\nabla_{V}^{N} \phi {F}_{\ast} X, U) = 0$ and \\$g_{N} (\nabla_{V}^{N} B V, U)$ $= 0$, for any $U \in \Gamma(range {F}_{\ast})^{\bot}, \nabla_{V}^{N} B V, \nabla_{V}^{N} \phi {F}_{\ast} X \in \Gamma(range {F}_{\ast})$. \\Also by using (\ref{2.2}), we have $\nabla_{{F}_{\ast}X}^{N} B V = \overset{N}{\nabla_{{F}_{\ast}X}^{{F}}} BV \circ {F} = (\nabla {F}_{\ast}) (X, {}^{\ast} {F}_{\ast} B V) + {F}_{\ast}(\nabla_{X}^{M \ast} {F}_{\ast} B V)$. In addition, since $(range {F}_{\ast})^{\bot}$ is totally geodesic, we have $\nabla_{V}^{N} \omega {F}_{\ast} X = \nabla_{V}^{{F} \bot} \omega {F}_{\ast}X$ and $\nabla_{V}^{N} C V = \nabla_{V}^{{F} \bot} C V$. Then by (\ref{3.8}), we obtain
		\begin{equation*}
			\nabla_{\dot{\beta}}^{N} \dot{\beta} = - P\left( \begin{array}{c} - S_{\omega {F}_{\ast}X} {F}_{\ast} X +{F}_\ast(\nabla_X^{M} {}^{\ast}{F}_\ast \phi {F}_\ast X) + (\nabla {F}_\ast) (X, {}^{\ast}{F}_\ast \phi {F}_\ast X) \\+ \nabla_{X}^{{F}\bot} \omega {F}_{\ast}X + (\nabla {F}_{\ast}) (X,{}^{\ast} {F}_{\ast} B V ) + {F}_{\ast}(\nabla_{X}^{M\ast} {F}_{\ast}BV) \\- S_{CV} {F}_{\ast} X  + \nabla_{X}^{{F} \bot} C V + \nabla_{V}^{N} \phi {F}_{\ast} X \\+ \nabla_{V}^{{F} \bot}\omega{F}_{\ast} X + \nabla_{V}^{N} B V + \nabla_{V}^{{F} \bot}C V
			\end{array}
			\right).
		\end{equation*}
		Now, $\beta$ is geodesic on $N \iff \nabla_{\dot{\beta}}^{N} \dot{\beta} = 0 \iff$
		\begin{eqnarray*}
			0 & = & - S_{\omega {F}_{\ast}X} {F}_{\ast} X +{F}_\ast(\nabla_X^{M} {}^{\ast}{F}_\ast \phi {F}_\ast X) + (\nabla {F}_\ast) (X, {}^{\ast}{F}_\ast \phi {F}_\ast X) \\&&+ \nabla_{X}^{{F}\bot}\omega {F}_{\ast} X + (\nabla {F}_{\ast}) (X, {}^{\ast} {F}_{\ast} B V) + {F}_{\ast}(\nabla_{X}^{M \ast} {F}_{\ast} B V) - S_{CV} {F}_{\ast} X \\&&+ \nabla_{X}^{{F} \bot} C V  + \nabla_{V}^{N} \phi {F}_{\ast} X + \nabla_{V}^{{F} \bot} \omega {F}_{\ast} X + \nabla_{V}^{N} B V + \nabla_{V}^{{F} \bot} C V.
		\end{eqnarray*}
		By separating the vertical and horizontal parts of above equation, we obtain (\ref{3.5}) and (\ref{3.6}). This completes the proof.
	\end{proof}
	
	\begin{theorem}
		Let ${F}$ be a SIRM from a Riemannian manifold $(M, g_{M})$ to a K\"ahler manifold $(N, g_{N}, P)$ such that leaves of $range{F}_\ast$ are connected, and $\alpha, \beta ={F} \circ \alpha$ be geodesic curves on $M$ and $N$, respectively. Then, ${F}$ is CSIRM with $s=e^{f}$ if and only if 
		\begin{align*}
			\left\Vert {F}_{\ast} X \right\Vert^{2} g_{N} (\nabla^{N}f, V)& = g_{N} (S_{\omega {F}_{\ast}X} {F}_{\ast} X + S_{CV}{F}_{\ast} X -\nabla_{V}^{N} \phi {F}_{\ast}X \\&- {F}_\ast(\nabla_X^{M} {}^{\ast}{F}_\ast \phi {F}_\ast X), BV) - g_{N}((\nabla {F}_\ast) (X, {}^{\ast}{F}_\ast \phi {F}_\ast X) \\&+ \nabla_{X}^{{F} \bot}\omega {F}_{\ast} X + (\nabla {F}_{\ast}) (X, {}^{\ast} {F}_{\ast} B V) + \nabla_{V}^{{F} \bot} \omega {F}_{\ast}X, C V),
		\end{align*}
		where $f$ is a smooth function on $N$ and ${F}_{\ast} X \in \Gamma(range {F}_{\ast}), V \in \Gamma(range {F}_{\ast})^{\bot}$ are components of $\dot{\beta}(t)$.
	\end{theorem}
	
	\begin{proof}
		Let $\alpha$ be a geodesic curve on $M$ and $\beta = {F} \circ \alpha$ be geodesic on $N$ with ${F}_{\ast} X \in \Gamma (range {F}_{\ast})$ and $V \in \Gamma(range {F}_{\ast})^{\bot}$ are components of $\dot{\beta}(t)$, and $\theta (t)$ denote the angle in $\left[0, \pi \right]$ between $\dot{\beta}$ and $V$. Let $\sqrt{k}$ be constant speed of $\beta$ on $N$ that is, $k = g_{N} (\dot{\beta}(t), \dot{\beta}(t)) = \left\Vert \dot{\beta}(t) \right\Vert^{2}$. Hence we conclude that,
		\begin{equation}\label{3.9}
			g_{N}({F}_{\ast} X, {F}_{\ast} X) = k \sin^{2} \theta (t)
		\end{equation}
		and
		\begin{equation}\label{3.10}
			g_{N}(V, V) = k \cos^{2} \theta (t).
		\end{equation}
		Differentiating (\ref{3.10}), we have
		\begin{equation}\label{3.11}
			\frac{d}{dt} g_{N} (V, V) = - 2 k \sin \theta (t) \cos \theta (t) \frac{d \theta}{dt}. 
		\end{equation}
		On the other hand by (\ref{2.8}), we have 
		\begin{equation}\label{3.12}
			\frac{d}{dt} g_{N} (V, V)=\frac{d}{dt} g_{N} (PV, PV). 
		\end{equation}
		Using (\ref{3.4}) in (\ref{3.12}), we get
		\begin{equation*}
			\frac{d}{dt} g_{N} (V, V) = \frac{d}{dt} \left(g_{N} (BV, BV) + g_{N} (CV, CV)\right),
		\end{equation*}
		which implies
		\begin{equation}\label{3.13}
			\frac{d}{dt} g_{N}(V, V) = 2 g_{N} (\nabla_{\dot{\beta}}^{N}BV, BV) + 2 g_{N}(\nabla_{\dot{\beta}}^{N} C V, C V).
		\end{equation}
		Putting $\dot{\beta} = {F}_{\ast} X + V$ in (\ref{3.13}), we get
		\begin{equation}\label{3.14}
			\begin{array}{ll}
				\frac{d}{dt} g_{N} (V, V) = & 2 g_{N} (\nabla_{{F}_{\ast} X}^{N} B V, B V) + 2 g_{N} (\nabla_{{F}_{\ast} X}^{N} C V,C V) \\& + 2 g_{N} (\nabla_{V}^{N} B V, B V) + 2 g_{N} (\nabla_{V}^{N} C V, C V).
			\end{array}
		\end{equation}
		Since $(range {F}_{\ast})^{\bot}$ is totally geodesic, (\ref{3.14}) can be written as 
		\begin{equation}\label{3.15}
			\begin{array}{ll}
				\frac{d}{dt} g_{N} (V, V) = & 2 g_{N} (\overset{N}{\nabla_{X}^{{F}}} BV \circ {F}, B V) + 2 g_{N}(\nabla_{{F}_{\ast} X}^{N}C V, C V) \\& + 2 g_{N} (\nabla_{V}^{N} B V, B V) + 2 g_{N}(\nabla_{V}^{{F} \bot} C V, C V).
			\end{array}
		\end{equation}
		Using (\ref{2.2}), (\ref{2.3}) and (\ref{2.5}) in (\ref{3.15}), we get
		\begin{align}\label{3.17}
			\frac{d}{dt} g_{N} (V, V) =& 2 g_{N} ({F}_{\ast}(\nabla_{X}^{M\ast} {F}_{\ast} B V) + \nabla_{V}^{N} B V, B V) \nonumber\\&+ 2 g_{N} (\nabla_{X}^{{F} \bot} C V + \nabla_{V}^{{F} \bot} C V, C V). 
		\end{align}
		Using (\ref{3.5}) and (\ref{3.6}) in (\ref{3.17}), we get
		\begin{align}\label{3.18}
			\frac{d}{dt} g_{N}(V, V) = & 2 g_{N} (S_{\omega {F}_{\ast} X } {F}_{\ast} X + S_{C V} {F}_{\ast} X - \nabla_{V}^{N} \phi {F}_{\ast} X \nonumber\\&- {F}_\ast(\nabla_X^{M} {}^{\ast}{F}_\ast \phi {F}_\ast X), B V) - 2 g_{N} ((\nabla {F}_\ast) (X, {}^{\ast}{F}_\ast \phi {F}_\ast X) \nonumber\\&+ \nabla_{X}^{{F} \bot} \omega {F}_{\ast} X + (\nabla {F}_{\ast})(X, {}^{\ast} {F}_{\ast} B V ) + \nabla_{V}^{{F} \bot} \omega {F}_{\ast} X, C V). 
		\end{align}
		From (\ref{3.11}) and (\ref{3.18}), we have
		\begin{align}\label{3.19}
			- k \sin \theta (t) \cos \theta (t) \frac{d \theta}{dt} = & g_{N} (S_{\omega {F}_{\ast}X} {F}_{\ast} X + S_{CV} {F}_{\ast} X - \nabla_{V}^{N} \phi {F}_{\ast}X \nonumber \\&- {F}_\ast(\nabla_X^{M} {}^{\ast}{F}_\ast \phi {F}_\ast X), BV) \nonumber \\&- g_{N} ((\nabla {F}_\ast) (X, {}^{\ast}{F}_\ast \phi {F}_\ast X) \nonumber \\&+ \nabla_{X}^{{F} \bot} \omega {F}_{\ast} X + (\nabla{F}_{\ast}) (X, {}^{\ast} {F}_{\ast} B V) \nonumber \\&+ \nabla_{V}^{{F} \bot} \omega {F}_{\ast} X, C V).
		\end{align}
		Now, ${F}$ is a CSIRM with $s = e^{f}$ if and only if $\frac{d}{dt} (e^{f\circ \beta} \sin \theta (t)) = 0$. Therefore
		\begin{equation*}
			\frac{d}{dt} (e^{f\circ \beta} \sin \theta (t)) = 0 \iff e^{f \circ \beta} \frac{d(f \circ \beta)}{dt} \sin \theta (t) + e^{f \circ \beta} \cos \theta (t) \frac{d \theta}{dt} = 0.
		\end{equation*}
		Since $s$ is a positive function,
		\begin{equation*}
			\frac{d(f\circ \beta)}{dt} \sin \theta (t) + \cos \theta (t) \frac{d \theta}{dt} = 0.
		\end{equation*}
		By multiplying this with non-zero factor $k \sin \theta$ and using (\ref{3.9}), we get
		\begin{equation}\label{3.20}
			-k \sin \theta (t) \cos \theta (t) \frac{d \theta }{dt} = g_{N}({F}_{\ast} X, {F}_{\ast} X) \frac{d(f \circ \beta)}{dt}. 
		\end{equation}
		Since the left-hand sides of (\ref{3.19}) and (\ref{3.20}) are equal, we can write
		\begin{align*}
			g_{N}({F}_{\ast}X, {F}_{\ast} X) g_{N}(\nabla^N f, V) & = g_{N} (S_{\omega{F}_{\ast} X} {F}_{\ast} X + S_{CV}{F}_{\ast} X - \nabla_{V}^{N} \phi {F}_{\ast} X \\&- {F}_\ast(\nabla_X^{M} {}^{\ast}{F}_\ast \phi {F}_\ast X), B V) \\&- g_{N}((\nabla {F}_\ast) (X, {}^{\ast}{F}_\ast \phi {F}_\ast X) \\&+ \nabla_{X}^{{F} \bot} \omega {F}_{\ast} X + (\nabla {F}_{\ast}) (X, {}^{\ast} {F}_{\ast} B V) \\&+ \nabla_{V}^{{F} \bot} \omega {F}_{\ast} X, C V).
		\end{align*}
		Since $\frac{d(f\circ \beta)}{dt} = g_{N} (\nabla^N f, \dot{\beta}(t)) = g_{N} (\nabla^N f, {F}_{\ast} X + V) = g_{N}(\nabla^N f, V)$. Hence the theorem is proved.
	\end{proof}
	
	\begin{proposition}
		Let ${F}$ be a CSIRM with $s=e^{f}$ from a Riemannian manifold $(M, g_{M})$ to a K\"ahler manifold $(N, g_{N}, P)$. Then, following statements are true:
		\begin{enumerate}[$(i)$]
			\item $\nabla^N_{P F_\ast X} F_\ast Y = -g_N(F_\ast Y, \phi F_\ast X) \nabla^N f + F_\ast(\nabla^M_{{}^{\ast}F_\ast \phi F_\ast X} Y)~\\ \text{if $F_\ast X, F_\ast Y \in \Gamma(D_1)$}$.
			
			\item $\nabla^N_{P F_\ast X} F_\ast Y = F_\ast(\nabla^M_{{}^{\ast}F_\ast \phi F_\ast X} Y)~\text{if $F_\ast X \in \Gamma(D_1)$ and $F_\ast Y \in \Gamma(D_2)$}$.
			
			\item $\nabla^N_{PU} V = V(f) BU + \nabla^{F \perp}_{{}^{\ast}F_\ast BU} V + \nabla^{F \perp}_{CU} V~\text{if $U, V \in \Gamma(rangeF_\ast)^\perp$}$.
			
			\item $\nabla^N_{U} PV = U(f) BU + \nabla^{F \perp}_{{}^{\ast}F_\ast BV} U + \nabla^{F \perp}_{U} CV~\text{if $U, V \in \Gamma(rangeF_\ast)^\perp$}$.
		\end{enumerate}
	\end{proposition}
	
	\begin{proof}
		For $F_\ast X, F_\ast Y \in \Gamma(D_1)$:
		\begin{align*}
			\nabla^N_{P F_\ast X} F_\ast Y &=\nabla^N_{\phi F_\ast X} F_\ast Y = \nabla^N_{{}^{\ast}F_\ast \phi F_\ast X} F_\ast Y\\& = (\nabla F_\ast)(Y, {}^{\ast}F_\ast \phi F_\ast X) + F_\ast(\nabla^M_{{}^{\ast}F_\ast \phi F_\ast X} Y) \\& = -g_M(Y, {}^{\ast}F_\ast \phi F_\ast X)\nabla^N f + F_\ast (\nabla^M_{{}^{\ast}F_\ast \phi F_\ast X} Y)\\& = -g_N(F_\ast Y, \phi F_\ast X)\nabla^N f + F_\ast (\nabla^M_{{}^{\ast}F_\ast \phi F_\ast X} Y).
		\end{align*}
		Similarly for $F_\ast X \in \Gamma(D_1)$ and $F_\ast Y \in \Gamma(D_2)$:
		\begin{align*}
			\nabla^N_{P F_\ast X} F_\ast Y & = -g_N(F_\ast Y, \phi F_\ast X)\nabla^N f + F_\ast (\nabla^M_{{}^{\ast}F_\ast \phi F_\ast X} Y) \\&= F_\ast (\nabla^M_{{}^{\ast}F_\ast \phi F_\ast X} Y).
		\end{align*}
		Now, for $U, V \in \Gamma(rangeF_\ast)^\perp$:
		\begin{align*}
			\nabla^N_{PU} V& = \nabla^N_{BU + CU} V = \nabla^N_{BU} V + \nabla^{F \perp}_{CU} V\\& = -S_V BU + \nabla^{F \perp}_{{}^{\ast}F_\ast BU} V + \nabla^{F \perp}_{CU} V\\& = V(f) BU + \nabla^{F \perp}_{{}^{\ast}F_\ast BU} V + \nabla^{F \perp}_{CU} V.
		\end{align*}
		Similarly, we can obtain
		\begin{align*}
			\nabla^N_{U} PV = U(f) BV + \nabla^{F \perp}_{{}^{\ast}F_\ast BV} U + \nabla^{F \perp}_{U} CV.
		\end{align*}
		This completes the required proof.
	\end{proof}
	
	\begin{theorem}\label{thmharmonic} \cite{Meena}
		Let $F$ be a CRM with $s = e^f$ between Riemannian manifolds $(M, g_M)$ and $(N, g_N)$ such that $kerF_\ast$ is minimal. Then $F$ is harmonic if and only if $f$ is constant.
	\end{theorem}
	
	\begin{theorem}
		Let ${F}$ be a CSIRM with $s=e^{f}$ from a Riemannian manifold $(M, g_{M})$ to a K\"ahler manifold $(N, g_{N}, P)$ such that $kerF_\ast$ is minimal. Then $F$ is harmonic if and only if
		\begin{equation*}
			\textrm{\normalfont trace} \left\{\omega S_{\omega F_{\ast }X}F_{\ast }X-C\nabla
			_{X}^{F\perp}\omega F_{\ast }X-(\nabla _{X}^{N}P\phi F_{\ast
			}X)^{(rangeF_{\ast })^{\perp }}\right\} = 0,
		\end{equation*}
		where $X\in \Gamma(kerF_{\ast})^{\perp}$.
	\end{theorem}
	
	\begin{proof}
		For $X\in \Gamma(kerF_{\ast})^{\perp}$ by using (\ref{2.2}),
		(\ref{3.3}), (\ref{2.5}) and (\ref{3.4}) we have
		\begin{eqnarray*}
			(\nabla F_{\ast })(X,X) &=&\overset{N}{\nabla_{X}^{{F} }} {F}_{\ast} X - {F}_{\ast}(\nabla_{X}^{M}X) \\
			&=&-P\overset{N}{\nabla_{X}^{{F}}} P{F}_{\ast} X - {F}_{\ast}(\nabla_{X}^{M}X) \\
			&=&-\nabla_{F_\ast X}^{N} P\phi F_{\ast }X+\phi S_{\omega F_{\ast }X}F_{\ast
			}X+\omega S_{\omega F_{\ast }X}F_{\ast }X \\
			&&-B\nabla_{X}^{F\perp }\omega F_{\ast }X-C\nabla _{X}^{F\perp }\omega
			F_{\ast }X-F_{\ast }(\nabla _{X}^{M}X).
		\end{eqnarray*}
		We know that the $rangeF_{\ast}$ part of $(\nabla F_{\ast })(X,X)$ is equal to $0$, i.e.
		\begin{equation*}
			\phi S_{\omega F_{\ast} X} F_{\ast} X - B \nabla_{X}^{F \perp} \omega F_{\ast}X - F_{\ast} (\nabla_{X}^{M}X) - (\nabla_{X}^{N} P \phi F_{\ast}X)^{rangeF_{\ast}} = 0.
		\end{equation*}
		Therefore,
		\begin{align}\label{3.39}
			&(\nabla F_{\ast })(X, X) = -g_M(X, X) \nabla^{N}f \nonumber \\&= \omega S_{\omega F_{\ast}X} F_{\ast}X - C \nabla_{X}^{F \perp} \omega F_{\ast} X - (\nabla_{X}^{N} P \phi
			F_{\ast}X)^{(rangeF_{\ast})^{\perp}}. 
		\end{align}
		Taking trace of (\ref{3.39}), we get
		{\footnotesize\begin{align*}
				\nabla^{N}f=-\frac{1}{(m-r)} \textrm{trace} \left\{\omega S_{\omega F_{\ast }X}F_{\ast }X-C\nabla
				_{X}^{F\perp}\omega F_{\ast }X-(\nabla _{X}^{N}P\phi F_{\ast
				}X)^{(rangeF_{\ast })^{\perp }}\right\}.
		\end{align*}}
		Since $F$ is a CRM with minimal fibers then proof follows by Theorem \ref{thmharmonic}.
	\end{proof}
	
	\begin{lemma}\label{Lemma3.1}
		Let ${F}$ be a CSIRM with $s=e^{f}$ from a Riemannian manifold $(M, g_{M})$ to a K\"ahler manifold $(N, g_{N}, P)$. Then, either $f$ is constant on $P(D_{2})$ or $D_2$ is one dimensional.
	\end{lemma}
	
	\begin{proof}
		Since ${F}$ is CRM with $s=e^{f}$ then using Theorem \ref{TH1} and (\ref{2.2}) in (\ref{2.4}), we have
		\begin{equation}\label{3.21}
			\overset{N} {\nabla_{X}^{{F}}} {F}_{\ast} Y - {F}_{\ast} (\nabla_{X}^{M} Y) = - g_{M}(X, Y) \nabla^{N} f
		\end{equation}
		for $X, Y\in \Gamma (\bar{D_2})$. Taking the inner product of (\ref{3.21}) with $P {F}_{\ast} Z \in \Gamma (PD_{2})$, we have
		\begin{equation}\label{3.22}
			g_{N} (\overset{N}{\nabla_{X}^{{F}}} {F}_{\ast} Y - {F}_{\ast}(\nabla_{X}^{M}Y), P{F}_{\ast} Z) = - g_{M}(X, Y)g_{N}(\nabla^{N}f, P{F}_{\ast} Z).
		\end{equation}
		Since $\overset{N}{\nabla^{{F}}}$ is pull-back connection of the Levi-Civita connection ${\nabla}^{N}$, $\overset{N}{\nabla^{{F}}}$ is also Levi-Civita connection. Then, using compatibility condition in (\ref{3.22}), we get
		\begin{equation}\label{3.23}
			g_{N} (P\overset{N}{\nabla_{X}^{{F}}} {F}_{\ast} Z, {F}_{\ast} Y) = g_{M}(X, Y) g_{N}(\nabla^{N}f, P {F}_{\ast} Z).
		\end{equation}
		Using (\ref{2.8}) in (\ref{3.23}), we get
		\begin{equation}\label{3.24}
			-g_{N} (\overset{N}{\nabla_{X}^{{F}}} {F}_{\ast} Z, P {F}_{\ast} Y) = g_{M} (X, Y) g_{N}(\nabla^{N}f, P {F}_{\ast} Z).
		\end{equation}
		Using (\ref{3.21}) in (\ref{3.24}), we get
		\begin{equation}\label{3.25}
			g_{M} (X, Z) g_{N} (\nabla^{N}f, P {F}_{\ast} Y) = g_{M} (X, Y) g_{N} (\nabla^{N}f, P {F}_{\ast} Z).
		\end{equation}
		Now, putting $X = Y$ in (\ref{3.25}), we get
		\begin{equation}\label{3.26}
			g_{M} (X, Z) g_{N} (\nabla^{N}f, P {F}_{\ast} X ) = g_{M} (X, X) g_{N} (\nabla^{N}f, P{F}_{\ast} Z).
		\end{equation}
		Now interchanging $X$ and $Z$ in (\ref{3.26}), we get
		\begin{equation}\label{3.27}
			g_{M} (X, Z) g_{N} (\nabla^{N}f, P {F}_{\ast} Z) = g_{M} (Z, Z) g_{N} (\nabla^{N}f, P {F}_{\ast} X).
		\end{equation}
		From (\ref{3.26}) and (\ref{3.27}), we have
		\begin{equation}\label{3.28}
			g_{N} (\nabla^{N}f, P {F}_{\ast}X) \left(1 - \frac{g_{M} (X, X) g_{M}(Z, Z)} {g_{M} (X, Z) g_{M} (X, Z)}\right) = 0.
		\end{equation}
		Since $\nabla^{N} f \in \Gamma (range {F}_{\ast})^\bot$, (\ref{3.28}) implies the proof.
	\end{proof}
	
	\begin{theorem}
		Let ${F}$ be a CSIRM with $s=e^{f}$ from a Riemannian manifold $(M, g_{M})$ to a K\"ahler manifold $(N, g_{N}, P)$ such that $\dim(D_2) > 1$. Then, ${F}$ is totally geodesic if and only if
		\begin{enumerate}[$(i)$]
			\item $ker {F}_{\ast}$ is totally geodesic, and
			
			\item $(ker {F}_{\ast})^\bot$ is totally geodesic, and
			
			\item for $X, Y, Z \in \Gamma(ker {F}_{\ast})^\bot$ such that ${F}_\ast Z = \phi {F}_\ast Y$, we have	
			\begin{align*}
				&g_{M} (X, Z)\left\{B(\nabla^{N} f) + C(\nabla^{N} f)\right\} = \phi {F}_\ast(\nabla_X^{M} Z) \\&+ \omega {F}_\ast(\nabla_X^{M} Z) + B \nabla_{X}^{{F} \bot} \omega {F}_\ast Y + C \nabla_{X}^{{F} \bot} \omega {F}_\ast Y + {F}_\ast(\nabla_X^{M} Y).
			\end{align*}
		\end{enumerate}
	\end{theorem}
	
	\begin{proof}
		We know that ${F}$ is totally geodesic if and only if 
		\begin{equation}\label{3.40}
			( \nabla {F}_{\ast})(U, V) = 0,
		\end{equation}
		\begin{equation}\label{3.41}
			( \nabla {F}_{\ast})(X, U) = 0
		\end{equation}
		and
		\begin{equation}\label{3.42}
			( \nabla {F}_{\ast})(X, Y) = 0
		\end{equation}	
		for $U, V \in \Gamma(ker{F}_{\ast})$ and $X, Y \in \Gamma(ker{F}_{\ast})^\bot$. More precisely, (\ref{3.40}) and (\ref{3.41}) have mean $ker{F}_{\ast}$ and $(ker{F}_{\ast})^\bot$ are totally geodesic. In addition, (\ref{3.42}) has means
		\[(\nabla {F}_\ast)(X, Y) = {\nabla_{X}^{{F} }} {F}_{\ast} (Y) - {F}_{\ast}(\nabla_{X}^{M}Y) = -P {\nabla_{X}^{{F}}} P {F}_{\ast} (Y) - {F}_{\ast}(\nabla_{X}^{M}Y) =0.\]
		Then, by (\ref{3.3}), we have
		\[-P\nabla_X^{{F}} (\phi {F}_{\ast} Y) - P \nabla_X^{{F}} \omega {F}_{\ast} Y - {F}_{\ast}(\nabla_X^{M} Y) =0.\]
		First using $\phi {F}_\ast Y = {F}_\ast Z$ and (\ref{2.5}) in above equation, and then using (\ref{2.2}), we obtain
		\begin{align*}
			&P((\nabla {F}_\ast)(X, Z) + {F}_{\ast} (\nabla_X^{M} Z)) \\&+ P(-S_{\omega {F}_{\ast} Y} {F}_{\ast} X + \nabla_X^{{F} \bot} \omega {F}_{\ast}Y) + {F}_\ast(\nabla_X^{M} Y) =0.
		\end{align*}
		Since ${F}$ is CRM, by using Theorem \ref{TH1} and Lemma \ref{lemmaforumbilicity}, we obtain
		\begin{align*}
			P(g_{M} (X, Z) \nabla^{N} f) &= P{F}_{\ast} (\nabla_X^{M} Z) + P( \omega {F}_{\ast} Y)(f) {F}_\ast X \\&+ P \nabla_X^{{F} \bot} \omega {F}_{\ast}Y + {F}_\ast(\nabla_X^{M} Y).
		\end{align*}
		Thus, proof follows by (\ref{3.3}) and Lemma \ref{Lemma3.1}.
	\end{proof}
	
	\begin{theorem}
		Let ${F}$ be a CSIRM with $s=e^{f}$ from a Riemannian manifold $(M, g_{M})$ to a K\"ahler manifold $(N, g_{N}, P)$ such that $\dim(D_2) > 1$. Then, following statements are true:
		\begin{enumerate}[$(i)$]
			\item $\nabla^{N} f \in \Gamma(P D_2)$.
			
			\item $S_\xi = 0$, for all $\xi \in \Gamma(\eta)$.
			
			\item $\nabla^N_{P F_\ast X} F_\ast Y = \nabla_X^{F \perp} P F_\ast Y$ \text{if $F_\ast X, F_\ast Y \in \Gamma(D_2)$}.

			\item $\nabla^N_{P F_\ast X} F_\ast Y = \nabla^N_{\phi F_\ast X} F_\ast Y + \nabla^{F \perp}_Y \omega F_\ast X$ \text{if $F_\ast X, F_\ast Y \in \Gamma(rangeF_\ast)$}.
		\end{enumerate}
	\end{theorem}
	
	\begin{proof}
		Since ${F}$ is CRM, for $X \in \Gamma(\bar{D_1})$, we have
		\begin{equation}\label{3.29}
			\nabla_X^{{F}} {F}_\ast X - {F}_{\ast}(\nabla_X^{M} X) = - g_{M}(X, X) \nabla^{N} f.
		\end{equation}
		Taking inner product with $V \in \Gamma(\eta)$, we obtain
		\begin{equation*}
			g_{N} (\nabla_X^{{F}} {F}_\ast X, V) = - g_{M}(X, X) g_{N}(\nabla^{N} f, V).
		\end{equation*}
		Since $N$ is K\"ahler manifold, we can write
		\begin{equation*}
			g_{N} (P \nabla_X^{{F}} {F}_\ast X, PV) = g_{N} (\nabla_X^{{F}} P {F}_\ast X, PV) = - g_{M}(X, X) g_{N}(\nabla^{N} f, V).
		\end{equation*}
		By using (\ref{3.3}) in above equation, we get
		\begin{equation}\label{3.30}
			g_{N} (\nabla_X^{{F}} \phi {F}_\ast X, PV) = - g_{M}(X, X) g_{N}(\nabla^{N} f, V).
		\end{equation}
		Using (\ref{3.29}) in (\ref{3.30}), we obtain
		\begin{equation*}
			g_{M}(X, {}^\ast {F}_{\ast}(\phi {F}_{\ast}X)) g_{N} (\nabla^{N}f, PV) = g_{M}(X, X) g_{N}(\nabla^{N} f, V).
		\end{equation*}
		Since ${F}$ is RM, we can write 
		\begin{equation}\label{3.31}
			g_{N}({F}_{\ast} X, \phi {F}_{\ast}X) g_{N} (\nabla^{N}f, PV) = g_{M}(X, X) g_{N}(\nabla^{N} f, V).
		\end{equation}
		Interchanging $V$ and $\nabla^{N} f$, we get
		\begin{equation}\label{3.32}
			g_{N}({F}_{\ast} X, \phi {F}_{\ast}X) g_{N} (V, P\nabla^{N}f) = g_{M}(X, X) g_{N}(\nabla^{N} f, V).
		\end{equation}
		Using (\ref{3.32}) in (\ref{3.31}), we obtain
		\begin{equation*}
			g_{N}({F}_{\ast} X, \phi {F}_{\ast}X) g_{N} (PV, \nabla^{N}f) = 0.
		\end{equation*}
		This implies the proof of $(i)$.
		For $\xi \in \Gamma(\eta)$ and $F_\ast X \in \Gamma(D_2)$:
		\begin{align*}
			g_N (S_\xi F_\ast X, F_\ast X) & = - \xi (f) g_N(F_\ast X, F_\ast X) \\& = -g_N(\xi, \nabla^N f) g_N(F_\ast X, F_\ast X).
		\end{align*}
		Since $\nabla^N f \in \Gamma(PD_2)$ and $(\eta) \perp (PD_2)$, by above equation we obtain
		\begin{align*}
			g_N (S_\xi F_\ast X, F_\ast X) = 0.
		\end{align*}
		Similarly for $\xi \in \Gamma(\eta)$ and $F_\ast Y \in \Gamma(D_1)$, we obtain
		\begin{align*}
			g_N (S_\xi F_\ast Y, F_\ast Y) = 0.
		\end{align*}
		Thus, $S_\xi F_\ast X = 0 = S_\xi F_\ast Y$ for $\xi \in \Gamma(\eta), F_\ast Y \in \Gamma(D_1)$ and $F_\ast X \in \Gamma(D_2)$. Hence $S_\xi = 0$ on $(range F_\ast)$. Thus the second statement hold.
		The next statements can easily prove by using (\ref{2.5}), $\nabla^{N} f \in \Gamma(P D_2)$ and Theorem \ref{TH1}.
	\end{proof}
	
	\begin{proposition}\label{prop1}
		Let ${F}$ be a CSIRM with $s=e^{f}$ from a Riemannian manifold $(M, g_{M})$ to a K\"ahler manifold $(N, g_{N}, P)$ such that $\dim(D_2) > 1$. Then the following statements are true:
		\begin{enumerate}[$(i)$]
			\item $\bar{D_1}$ defines a totally geodesic foliation if and only if $(\nabla {F}_{\ast})(X, U)$ has no component in $D_1$ for $X \in \Gamma(\bar{D_1})$ and $U \in \Gamma(ker {F}_\ast)$.
			
			\item $\bar{D_2}$ defines a totally geodesic foliation if and only if $(\nabla {F}_{\ast})(X_2, U)$ has no component in $D_2$ for $X_2 \in \Gamma(\bar{D_2})$ and $U \in \Gamma(ker {F}_\ast)$.
		\end{enumerate}
	\end{proposition}
	\begin{proof}
		We know that $\bar{D_1}$ defines a totally geodesic foliation if and only if $g_{M}(\nabla_X^{M} Y, U) = 0$ and $g_{M} (\nabla_X^{M} Y, X') = 0$ for $X, Y \in \Gamma(\bar{D_1})$, $U \in \Gamma(ker {F}_\ast)$ and $X' \in \Gamma(\bar{D_2})$. Now since ${F}$ is RM, using (\ref{2.1}) and (\ref{2.2}) we can write
		\begin{align*}
			g_{M}(\nabla_X^{M} Y, U) = - g_{M} (\nabla_X^{M} U, Y) &= -g_{N}({F}_\ast (\nabla_X^{M} U), {F}_\ast Y) \\&= g_{N}((\nabla {F}_\ast)(X, U), {F}_\ast Y),
		\end{align*}
		and similarly 
		\begin{align*}
			g_{M}(\nabla_X^{M} Y, X')& = - g_{M} (\nabla_X^{M} X', Y) = -g_{N}({F}_\ast (\nabla_X^{M} X'), {F}_\ast Y) \\&= -g_{N}(\nabla_X^{{F}}{F}_\ast X', {F}_\ast Y)=-g_{N}(\nabla_{{F}_{\ast} X}^{N}{F}_\ast X', {F}_\ast Y).
		\end{align*}
		Since $N$ is K\"ahler manifold, by using (\ref{2.8}) and then (\ref{2.5}), we obtain
		\begin{align}\label{3.33}
			g_{M}(\nabla_X^{M} Y, X') =-g_{N}(\nabla_{{F}_{\ast} X}^{N} P{F}_\ast X', P{F}_\ast Y) \nonumber\\=g_{N}(S_{P {F}_\ast X'} {F}_\ast X, P {F}_\ast Y).
		\end{align}
		Using Theorem \ref{TH1} and then Lemma \ref{Lemma3.1} in (\ref{3.33}), we get
		\begin{align*}
			g_{M}(\nabla_X^{M} Y, X') = -(P {F}_{\ast} X'(f))g_{N} ({F}_\ast X, P {F}_\ast Y) = 0.
		\end{align*}
		This implies the proof of $(i)$.
		
		On the other hand, we know that $\bar{D_2}$ defines totally geodesic if and only if $g_{M}(\nabla_{X_2}^{M} {Y_2}, U) = 0$ and $g_{M} (\nabla_{X_2}^{M} {Y_2}, {X_2}') = 0$ for ${X_2}, {Y_2} \in \Gamma(\bar{D_2})$, $U \in \Gamma(ker {F}_\ast)$ and ${X_2}' \in \Gamma(\bar{D_1})$. Now since ${F}$ is RM, using (\ref{2.1}) and (\ref{2.2}) we can write
		\begin{align*}
			g_{M}(\nabla_{X_2}^{M} {Y_2}, U) = - g_{M} (\nabla_{X_2}^{M} U, {Y_2}) &= -g_{N}({F}_\ast (\nabla_{X_2}^{M} U), {F}_\ast {Y_2}) \\&= g_{N}((\nabla {F}_\ast)({X_2}, U), {F}_\ast {Y_2}),
		\end{align*}
		and similarly 
		\begin{align*}
			g_{M}(\nabla_{X_2}^{M} {Y_2}, {X_2}')& = g_{N}({F}_\ast (\nabla_{X_2}^{M} {Y_2}), {F}_\ast {X_2}') \\&= g_{N}(\nabla_{X_2}^{{F}}{F}_\ast {Y_2}, {F}_\ast {X_2}')=g_{N}(\nabla_{{F}_{\ast} {X_2}}^{N}{F}_\ast {Y_2}, {F}_\ast {X_2}').
		\end{align*}
		Since $N$ is K\"ahler manifold, by using (\ref{2.8}) and then (\ref{2.5}), we obtain
		\begin{align}\label{3.34}
			g_{M}(\nabla_{X_2}^{M} {Y_2}, {X_2}') &= g_{N}(\nabla_{{F}_{\ast} {X_2}}^{N} P{F}_\ast {Y_2}, P{F}_\ast {X_2}') \nonumber\\&=-g_{N}(S_{P {F}_\ast {Y_2}} {F}_\ast {X_2}, P {F}_\ast {X_2}').
		\end{align}
		Using Theorem \ref{TH1} and then Lemma \ref{Lemma3.1} in (\ref{3.34}), we get
		\begin{align*}
			g_{M}(\nabla_{X_2}^{M} {Y_2}, {X_2}') = -(P {F}_{\ast} {Y_2}(f))g_{N} ({F}_\ast {X_2}, P {F}_\ast {X_2}') = 0.
		\end{align*}
		This implies the proof of $(ii)$.
	\end{proof}
	
	\begin{definition}\label{locallyproduct}
		\cite{Ponge_1993} Let $(M, g_M)$ be a Riemannian manifold and assume that the canonical foliations $L_1$ and $L_2$ such that $L_1 \cap L_2 = \{0\}$ everywhere. Then, $(M, g_M)$ is a locally product manifold if and only if $L_1$ and $L_2$ are totally geodesic foliations. 
	\end{definition}
	\begin{theorem}\label{thm3.5}
		Let ${F}$ be a CSIRM with $s=e^{f}$ from a Riemannian manifold $(M, g_{M})$ to a K\"ahler manifold $(N, g_{N}, P)$ such that $\dim(D_2) > 1$. Then $(ker{F}_{\ast})^\bot$ is a locally product manifold of $\bar{D_1}$ and $\bar{D_2}$ if and only if
		\begin{enumerate}[$(i)$]
			\item $(\nabla {F}_{\ast})(X, U)$ has no component in $D_1$, and
			
			\item $(\nabla {F}_{\ast})(X_2, U)$ has no component in $D_2$,
		\end{enumerate}
		for $X \in \Gamma(\bar{D_1}), X_2 \in \Gamma(\bar{D_2})$ and $U \in \Gamma(ker {F}_\ast)$.
	\end{theorem}
	
	\begin{proof}
		The proof is straight forward by Proposition \ref{prop1} and Definition \ref{locallyproduct}.
	\end{proof}
	
	\begin{theorem}\label{thm3.6}
		Let ${F}$ be a CSIRM with $s=e^{f}$ from a Riemannian manifold $(M, g_{M})$ to a K\"ahler manifold $(N, g_{N}, P)$ such that $\dim(D_2) > 1$. Then the base manifold is locally a product manifold ${{(range{F}_\ast)}} \times {{(range{F}_\ast)^\bot}}$ if and only if $g_{N} ({F}_{\ast}(\nabla_{X}^{M} \phi {F}_{\ast} X), BV) + g_{N} (\nabla_{X}^{{F} \bot} \omega {F}_{\ast} X, CV) + g_{N} (S_{CV} {F}_{\ast} X, \phi {F}_{\ast} X)= 0$ for $X \in \Gamma(\bar{D_1})$ and $V \in \Gamma(range {F}_\ast)^\bot$.
	\end{theorem}
	
	\begin{proof}
		Since $N$ is K\"ahler manifold, by using (\ref{2.8}), we have
		\begin{equation}\label{3.43}
			g_{N} (\nabla_{{F}_\ast X}^{N} {F}_\ast X, V) = g_{N}( \nabla_{{F}_\ast X}^{N} P{F}_\ast X, PV)
		\end{equation}
		for ${F}_{\ast} X \in \Gamma(range {F}_{\ast})$ and $V \in \Gamma(range {F}_{\ast})^\bot$.
		By using (\ref{3.3}) and (\ref{3.4}) in (\ref{3.43}), we write
		\begin{align*}
			g_{N} (\nabla_{{F}_\ast X}^{N} {F}_\ast X, V)&= g_{N}(\nabla_{{F}_\ast X}^{N} \phi {F}_\ast X, BV) + g_{N}(\nabla_{{F}_\ast X}^{N} \omega {F}_\ast X, BV) \\&+ g_{N}(\nabla_{{F}_\ast X}^{N} \phi {F}_\ast X, CV) + g_{N}(\nabla_{{F}_\ast X}^{N} \omega {F}_\ast X, CV).
		\end{align*}
		Using (\ref{2.5}) in above equation, we write
		\begin{align}\label{3.44}
			g_{N} (\nabla_{{F}_\ast X}^{N} {F}_\ast X, V)&= g_{N}(\nabla_{{F}_{\ast} X}^{N} \phi {F}_{\ast} X, BV) + g_{N} (S_{CV} {F}_{\ast} X, \phi {F}_{\ast}X) \nonumber \\&- g_{N}(S_{\omega {F}_{\ast} X} {F}_{\ast} X, BV) + g_{N}(\nabla_X^{{F} \bot} \omega {F}_{\ast} X, CV).
		\end{align}
		Using Theorem \ref{TH1} and then Lemma \ref{Lemma3.1} in (\ref{3.44}), we get
		\begin{align*}
			g_{N} (\nabla_{{F}_\ast X}^{N} {F}_\ast X, V)&= g_{N}(\nabla_{{F}_{\ast} X}^{N} \phi {F}_{\ast} X, BV) + g_{N}(\nabla_X^{{F} \bot} \omega {F}_{\ast} X, CV) \\&+ g_{N} (S_{CV} {F}_{\ast} X, \phi {F}_{\ast}X).
		\end{align*}
		Then, by using (\ref{2.2}), we get
		\begin{align*}
			g_{N} (\nabla_{{F}_\ast X}^{N} {F}_\ast X, V)&= g_{N}((\nabla {F}_{\ast})(X, {}^{\ast}{F}_{\ast}(\phi {F}_{\ast}X)), BV) \\&+ g_{N}(\nabla_X^{{F} \bot} \omega {F}_{\ast} X, CV) + g_{N} (S_{CV} {F}_{\ast} X, \phi {F}_{\ast}X).
		\end{align*}
		Then by Definition \ref{locallyproduct} proof is completed.
	\end{proof}
	
	\section{Further Studies with Almost Hermitian or K\"ahler Manifolds}\label{sec4}
	The Clairaut condition for RMs introduced in \cite{Sahin5} is investigated for IRMs, AIRMs and SIRMs in \cite{Yadav}, \cite{Yadav_mjm} and \cite{kpk_tjm}, respectively. On the other hand, the Clairaut condition for RMs introduced in \cite{Meena} is investigated for IRMs, AIRMs, SIRMs and SRMs \cite{Sahin_slant1} in \cite{Yadav}, \cite{Meena}, present paper and \cite{Jyoti} respectively. As further studies the Clairaut conditions can be investigated for (P)SRMs \cite{Sahin_slant2, AGslant, GAslant}, (P)SSRMs \cite{Park_semislant, ParkSahin, GAsemislant} and (P)HSRMs \cite{PHSRM-AG, Sahin_hemislant, GAS}.\\
	
	
	\noindent \textbf{Acknowledgments}\\
	Kiran Meena gratefully acknowledges the research facilities provided by Harish-Chandra Research Institute, Prayagraj, India. In addition, both authors are grateful to the referees for their comments to improve the paper.\\
	
	\noindent \textbf{Author contributions}
	Conceptualization, methodology, investigation, writing --- original and revised versions, validation, review, editing and reading have been performed by both the authors and agreed to the paper.\\
	
	\noindent \textbf{Funding}
	Kiran Meena is partial financial supported by the Department of Atomic Energy, Government of India [Offer Letter No.: HRI/133/1436 Dated 29 November 2022].\\
	
	\noindent \textbf{Data Availability}
	This manuscript has no associated data.\\
	
	\noindent \textbf{Declarations}\\\\
	\textbf{Conflict of Interest}
	The authors have no conflict of interest and no financial interests for this article.\\
	
	\noindent \textbf{Ethical Approval}
	The submitted work is original and not submitted to more than one journal for simultaneous consideration.\\
	
	\noindent \textbf{Publisher's Note}
	Springer Nature remains neutral with regard to jurisdictional claims in published maps and institutional affiliations.\\
	
	\noindent Springer Nature or its licensor (e.g. a society or other partner) holds exclusive rights to this article under a publishing agreement with the author(s) or other	rightsholder(s); author self-archiving of the accepted manuscript version of	this article is solely governed by the terms of such publishing agreement and	applicable law.\\

\end{document}